\documentclass[final,5p,times,twocolumn]{elsarticle}
\usepackage{amssymb}
\usepackage{amsfonts}
\usepackage{latexsym, amssymb,amsmath, amsbsy,amsopn, amstext}
\usepackage{amsmath, amsbsy}
\usepackage{hyperref}
\usepackage{float}
\usepackage{ifpdf}
\usepackage{xcolor}
\usepackage{tikz}
\usetikzlibrary{arrows,shapes,snakes,shadows,positioning,automata,patterns}
\usetikzlibrary{trees,decorations.pathmorphing,decorations.markings}
\def\0{{\bf 0}}

\newtheorem{thm}{Theorem}[section]
\newtheorem{lem}[thm]{Lemma}

\newtheorem{exa}[thm]{Example}
\newtheorem{rem}[thm]{Remark}

\newtheorem{assum}{Assumption}

\journal{Systems \& Control letters}

\begin{document}
\begin{frontmatter}

\title{    Distributed  Resource Allocation Over  Random  Networks  Based on Stochastic Approximation  \tnoteref{label0}}
\author{Peng Yi, Jinlong Lei, Yiguang Hong \corref{cor1}}
\ead{yipeng@amss.ac.cn, leijinlong11@mails.ucas.ac.cn, yghong@iss.ac.cn}
\address{The Key Laboratory  of Systems and Control, Academy of
Mathematics and Systems Science, Chinese Academy of Sciences}
\cortext[cor1]{Correspondence author: Yiguang Hong }

\begin{abstract}
In this paper, a stochastic approximation (SA)  based distributed algorithm is proposed to solve the resource allocation (RA)  with uncertainties. In this problem, a group of agents cooperatively  optimize a separable optimization problem with  a linear network resource constraint and {  allocation feasibility constraints},
where  the global objective function is the sum of agents' local objective functions.
{
Each agent can only get   noisy  observations of   its local function's gradient  and its local resource, which cannot be shared by other agents or transmitted to a center.  Moreover,  there are communication uncertainties such as time-varying topologies (described by random graphs) and additive channel noises.}
To solve the RA, we propose an SA-based distributed algorithm, and   prove that
 agents can collaboratively achieve the optimal allocation with probability one   by virtue of ordinary differential equation (ODE) method for SA.  Finally,
 simulations related to the  demand response management in power systems verify   the  effectiveness of  the proposed algorithm.
\end{abstract}
\tnotetext[label0]{This work was supported by Beijing Natural Science Foundation
(4152057), NSFC (61333001), and Program 973 (2014CB845301/2).}
\end{frontmatter}

\section{Introduction}

Resource allocation (RA) problem  is to  allocate the network resource   among  a group of agents  while optimizing
 certain performance index.
It has  drawn much research attention in many areas, such as the media access control in communication networks \cite{RAC}, signal processing in \cite{RAS}, and load demand management in \cite{stevon2}.
Hence, various RA models and RA algorithms have been proposed
(see \cite{RAC}-\cite{RA10} and the references therein).
However, most of existing algorithms need a center to collect the data over  networks or to coordinate  computation processes among all agents.

In fact, the center-free distributed optimization algorithms have attracted  more and more research attention in recent years \cite{Ned2}-\cite{peng}.
In  various network optimization problems,  the optimal decisions are made based on the whole network data, which,
however, are collected and stored by each individual agent of the network.
The distributed optimization algorithm keeps the data distributed through the network when seeking the optimal decision, and hence
eliminates the  ``one-to-all" communication burden and protects  agents' privacy.
Distributed optimization also endows each individual agent with autonomy and reactivity by  allowing it to formulate its local objective function and constraints with its local data.
From the network viewpoint,  the robustness to single  point failure  and the network scalability can be enhanced with distributed design.
Following the seminal work \cite{RA3}  of RA in large-scale networks along with the distributed optimization work in \cite{Ned2}-\cite{peng},
various center-free distributed algorithms for RA
have been proposed recently in \cite{DRA4}-\cite{DRA8}.

Stochastic approximation (SA) has been
adopted in distributed optimization algorithms
to address various kinds of uncertainties or to improve the computation efficiency.
In \cite{Ned3}, an SA-based distributed algorithm was
proposed when each agent can only get the noisy  observations of  its local  gradient, which extended the traditional SA optimization methods (see \cite{nemi}) to  distributed settings.
In \cite{DSA2},  an SA algorithm was given for distributed root seeking problem under  noisy observations, which was also a generalization of distributed optimization problems.
In practice, noisy gradient observations also exist in the zero-order distributed optimization algorithm as in \cite{Yuan}, and
randomized data sample was considered to reduce the computational  complexity in optimization with ``big data", resorting  to SA  for theoretical analysis (see \cite{DSA3}).  Besides,  SA algorithms were  also adopted  for distributed optimization to handle uncertainties in  communication systems in \cite{asum,ned}, and \cite{zhang2}.

Nevertheless, the existing distributed works of RA in \cite{DRA4}-\cite{DRA8}  have not considered various stochastic uncertainties related to information sharing  or data observations.
Since the problem  data is distributed throughout the network, each agent needs to share its local information with other agents through a communication network, which may involve various of  uncertainties.
Firstly, the communication network may switch  due to packet loss, media access control, or energy constraint.
To describe uncertainties of communication topologies,  different from the deterministic switching graphs in \cite{Ned2,lou12} and \cite{peng}, we adopt random graph models like \cite{asum,ned,boyd} and \cite{Zhang} here.
Secondly, the information shared through
the  network may not be accurate or may be corrupted by random noises due to  quantization errors or channel fading ( referring to \cite{peng}, \cite{zhang2} and \cite{Zhang}). On the other hand, noises can also  be actively added to the shared information for privacy protection as discussed  in \cite{DCU1}.
Moreover, agents may not get the exact local gradient or resource information due to  measurement or observation noises.

Main contributions of the paper are summarized as follows. (i)  A novel center-free distributed algorithm is  proposed to handle the RA problem, where each agent only utilizes noisy observations of  its local  gradient and resource information, and  noisy neighboring information shared through the randomly switching   networks.  (ii)
 The estimates are shown to converge  to the optimal allocation with probability one based on the ODE method for SA  algorithm.
(iii) The proposed model and algorithm are applied  to  distributed  multi-periods demand response management in power systems, along with  simulations  to show the effectiveness.

The remainder of the paper is organized as follows.
The RA  problem is formulated and an SA-based  distributed  algorithm is proposed in Section 2.
Then the  convergence result for the distributed algorithm is established in Section 3, while  simulation studies are shown in Section 4. Finally, the concluding remarks are given in Section 5.

\section{Problem Formulation and Proposed Algorithm}

Firstly, we  show related notations and preliminaries about convex analysis.
Denote $\mathbf{1}_m=(1,...,1)^T \in \mathbf{R}^m$ and
$\mathbf{0}_m=(0,...,0)^T \in \mathbf{R}^m$.
$col \{ x_1,\cdots, x_n\}= (x_1^T, \cdots, x_n^T)^T$ stacks  the vectors $x_1,  \cdots, x_n$.
$I_n$ denotes the identity matrix in $\mathbf{R}^{n\times n}$.
For a matrix $A=[a_{ij}]$, $a_{ij}$ or $A_{ij}$
stands for the matrix entry in the $i$th row and $j$th column of $A$.
$  \otimes $ denotes the  Kronecker  product. {  Denote $ker\{A\}$ and $range\{A\}$ as the null space and range space of matrix $A$, respectively.

For a  nonempty closed convex   set $\Omega \subset \mathbf{R}^m$ and a point $x \in \mathbf{R}^m$, denote $P_{\Omega} (x)$ as the point in $\Omega$ that is closest  to $x$, and call it the projection of $x$ on $\Omega$.
   $P_{\Omega} ( x)$ contains only one element for any $x \in  \mathbf{R}^m,$ and satisfies
       \begin{equation}\label{pro}
\| P_{\Omega} (x)-P_{\Omega} (y) \| \leq \| x-y\| ~~ \forall x, y\in\mathbf{R}^m.
\end{equation}
For a convex set $\Omega \subset \mathbf{R}^m$ and a point $x \in  \Omega$, define
   the normal cone  to $\Omega$ at $x$ as   $N_{\Omega} (x) \triangleq \{ v \in \mathbf{R}^m: \langle v, y-x \rangle \leq 0  ~~\forall y \in \Omega\}$. }

In the following two subsections, we formulate the distributed  RA  problem
with the  data observation and communication network models, and propose an SA-based distributed algorithm.
\subsection{Problem Formulation}
Consider a group of agents $\mathcal{N}=\{1, \cdots, n\}$ that cooperatively decide the optimal network resource allocation (RA), formulated  as follows:
{
\begin{equation}\label{problem}
\begin{split}
&\min_{x_i \in \mathbf{R}^m,i\in \mathcal{N}}  \qquad \sum_{i\in \mathcal{N}} f_i(x_i), \\
&subject \; to \;      \sum_{i\in \mathcal{N}} x_i = \sum_{i\in\mathcal{ N}} d_i, \quad x_i \in \Omega_i, i \in \mathcal{N}
\end{split}
\end{equation}
}
The   local allocation variable $x_i \in \mathbf{R}^m$ is decided by agent $i$, which is also associated with a local objective function $f_i(x_i)$. $d_i$ is the  local resource data, and can only be observed by agent $i$. The  resource of the whole network  is  the sum of all local resources, i.e.,  $\sum_{i \in \mathcal{N}}d_i$.
{ $\Omega_i$ is the local allocation feasibility constraint of agent $i$, and cannot be known by other  agents.  Furthermore,   $\Omega_i$ is  determined by $p_i$
inequality constraints:
$ \Omega_i=\{ x\in \mathbf{R}^m: q_{ij} (x) \leq 0, ~ \forall  j=1,\cdots, p_i\},$
where   $q_{ij}(\cdot),~j=1,\cdots, p_i$ are continuously differentiable   convex functions on
$\mathbf{R}^m$ .
Therefore, RA problem \eqref{problem} is  to find an allocation  that   minimizes the sum of local objective functions while satisfying the network resource constraint and the allocation feasibility constraints.}
The following assumptions can also be found in   \cite{RAC}-\cite{RA10}.

\begin{assum}\label{assp}
 Problem \eqref{problem} has a finite optimal solution.
 For any $i \in \mathcal{N}$, $f_i(x_i)$ is {    differentiable strictly convex function}, and moreover, its gradient  is globally  Lipschitz continuous, i.e.,    there exists a constant $l_c>0$ such that
$\| \nabla f_i(x) -\nabla f_i(y) \| \leq l_c \| x-y\|, \forall x,y \in \mathbf{R}^m .$
\end{assum}

{
The following constraint qualification assumption can be found in \cite{sa}.
\begin{assum}\label{assset}
For any $i \in \mathcal{N}, $ the set $\Omega_i$ is closed convex set and has nonempty interior points, and $\{ \nabla q_{ij}(x), ~ j \in \mathcal{I}_i(x)\}$ is linearly independent,  where $\mathcal{I}_i(x)=\{ j: q_{ij}(x)=0 \}$.
\end{assum}
}

The data observation model for agent $i$  at time $k$ is given as follows:
agent $i$ can get the   noisy  observation of its gradient $\nabla f_i(x_i)$ at given testing point $x_i(k)$  corrupted with noise $\nu_i(k)$ (that is,
$\nabla f_i(x_{i}(k)) + \nu_{i}(k)$) and the noisy local resource information corrupted with noise $\delta_i(k)$ (that is, $d_i+ \delta_i(k)$).
The stochastic gradient model should be  taken into consideration in  the following three cases:

(i) Stochastic optimization:
Agent  $i$'s local objective function takes the expectation form   as {  $f_i(x_i)= E_{\phi_i}[g(x_i, \phi_i)] = \int_{\Phi_i} g(x_i, \phi_i)d\mathbb{P}(\phi_i)$,}
where $\phi_i$ is a random vector supported on set $\Phi_i\in \mathbf{R}^d$  with probability distribution $\mathbb{P}$,  and $g_i : \mathbf{R}^m \times \Phi_i \rightarrow \mathbf{R}$.
%is well defined and finite valued for each  $x_i \in \mathbf{R}^m$.
It is more practical to utilize noisy gradient $\nabla g_i(x_i,\phi_i)$ given sampling $\phi_i$ rather than exact gradient  by performing multi-value integral at each iteration.
  In fact, the SA algorithm in \cite{nemi} and  DSA algorithm  in \cite{Ned3} considered this kind of gradient noise.

(ii) Zero-order optimization: When  agent $i$ can only get the  value of $f_i(x_i)$ given the testing point $x_i(k)$, the  gradient estimation methods, such as the Kiefer-Wolfowitz method in \cite{gfree1} and the randomized coordinate estimation in \cite{Yuan},  can lead to noisy gradient observations.

 (iii) Randomized data sample: If the local objective functions are constructed with ``big data",
a noisy gradient
based on randomly sampled data is  an  alternative to the exact gradient,
which may reduce the overall iteration computational complexity (see \cite{DSA3}).

Given the local data observations,  it is important and practical to solve \eqref{problem} in a distributed way,
where the agents need to share the local information with neighbors through switching  networks  and  noisy channels.

As we know, switching communication networks  can be  modeled by random graphs, e.g., \cite{asum}, \cite{ned}.
Denote a realization of the random graph at
time $k$ as $\mathcal{G}(k)=(\mathcal{N},\mathcal{E}(k))$, where
$\mathcal{E}(k) \subset \mathcal{N}\times \mathcal{N} $  is the  edge set at time $k$.
If agent $i$ can get  information from agent $j$ at time $k$, then $(j,i) \in \mathcal{E}(k)$ and
agent $j$ belongs to  agent  $i$'s neighbor set  $\mathcal{N}_i(k)=\{j|(j,i) \in
\mathcal{E}(k)\}$ at time $k$.
Define adjacency matrix $A(k)=[a_{ij}(k)]$ of $\mathcal{G}(k)$ with $a_{ij}(k)=1$ if
$j\in \mathcal{N}_i(k)$,  and $a_{ij}(k)=0$ otherwise.
Denote by
$Deg(k)=diag\{ \sum_{j=1}^n a_{1j}(k),..., \sum_{j=1}^n a_{nj}(k)\} $   the degree matrix, and by $L(k)=Deg(k)-A(k) $  the Laplacian  matrix of  $\mathcal{G}(k)$.

The following  assumption is given for  the random  graphs $\{\mathcal{G}(k)\}_{k \geq 1}$ (referring to \cite{asum}).

\begin{assum}\label{assg}
$\{L(k)\}$ is an i.i.d.  sequence with mean denoted by $\bar{L}= E[L(k)]$. Besides,  $\bar{L}$ is symmetric with $s_2(\bar{L})>0$,
where $s_2(\bar{L})$ denotes the secondly smallest eigenvalue of $\bar{L}$.
\end{assum}

% About the communication graph
\begin{rem}
Note that  Assumption \ref{assg}  does not require the communication graph to be connected or undirected at any time instance.
Only the mean graph is required to be  undirected and connected,
which ensures that  the local information can reach any other agents   in the average sense.
{  The gossip model in \cite{boyd} and  the broadcast model  in \cite{ned} are  also consistent with  Assumption \ref{assg}.}
\end{rem}

\subsection{SA-based Distributed Algorithm}

It is time to propose an SA-based distributed algorithm, based on assumptions  on  data observations  and communication  noises.

Denote   $x_i(k)$  as agent $i$'s estimate for its local optimal  allocation at time $k$,
and   denote $\lambda_i(k),~z_i(k)$  as the auxiliary variables of agent $i$.
 The agents share their  auxiliary variables through the communication network at each iteration.
If $(j,i)\in \mathcal{E}(k)$, then agent $i$ can get the noisy information of $\{\lambda_j(k),z_j(k)\}$, corrupted with noise
$ \zeta_{ij}(k)$ and $ \epsilon_{ij}(k)$,   from agent $j$.
Namely, $ \lambda_j(k)+\zeta_{ij}(k)$ and $z_j(k) +\epsilon_{ij}(k) $ are the values received by  agent $i$ from  agent $j$ at time $k$,   which are not separable.
Moreover, agent $i$ also has the local noisy gradient observation $\nabla f_i(x)+ \nu_{i}(k) $ and noisy resource observation $d_i+\delta_i(k)$.

The SA-based distributed recursive algorithm for agent $i$ is given as follows:
\begin{equation}\label{dy1}
\begin{array}{l} \hline
{\bf SA-based \; Distributed \; Resource \; Allocation \;Algorithm }        \\ \hline
\displaystyle {x}_{i}( k+1)            = P_{\Omega_i} \big( x_i(k)             + \alpha_{k} \big( -\big(\nabla f_i(x_i(k))+\nu_i(k)\big) + \lambda_{i }(k) \big) \big),\\
\displaystyle {\lambda}_{i}( k+1)      =   \lambda_{i}(k)    + \alpha_{k} \big( (d_{i } +\delta_{i}(k))- x_{i}(k) \\
\displaystyle  \qquad \qquad \quad        -  \sum_{j=1}^n a_{ij} (k)(\lambda_{i}(k)-(\lambda_{j }(k)+\zeta_{ ij }(k))) \\
\displaystyle  \qquad \qquad \quad        -  \sum_{j=1}^n a_{ij} (k)(z_{i}(k)-(z_{j}(k)+\epsilon_{ ij }(k)) ) \big), \\
\displaystyle{z}_i (k+1)              =  z_{i}(k)           + \alpha_{k}\sum_{j=1}^n a_{ij} (k)\big(\lambda_{i}(k)-(\lambda_{j}(k)+ \zeta_{ ij }(k)) \big),\\ \hline
\end{array}
\end{equation}
 where   the  step-size  $\{ \alpha_k\}$   satisfies
\begin{equation}\label{stepsize}
\alpha_k >0, ~~ \sum_{k=1}^{\infty}  \alpha_k =\infty, ~~\sum_{k=1}^{\infty}  \alpha_k^2 < \infty .
\end{equation}

Obviously, the algorithm \eqref{dy1} is a {\bf fully distributed} one  since each agent only
uses its local noisy observations and the noisy information received from  its neighbors, and only performs local projection with its local set $\Omega_i$.

{  Since the local objective functions $f_i(x_i)$  is convex and continuously differentiable,
the KKT condition of (\ref{problem}) is
\begin{equation}\label{kkt}
\begin{array}{lll}
&&\mathbf{0}_{m} \in   \nabla f_i(x_i^*)  -  \lambda^*  +  N_{\Omega_i}(x_i^*), i=1,\cdots,n \\
&& \sum_{i\in \mathcal{N}} x^*_i = \sum_{i\in\mathcal{ N}} d_i\; \quad x_i^* \in \Omega_i,
\end{array}
\end{equation}
Algorithm \eqref{dy1} is developed by combining  the ODE methods for KKT condition \eqref{kkt} and the ODE methods for stochastic approximation. In some sense,   $\lambda_i$ in \eqref{dy1} is the local ``copy" of  Lagrangian multiplier for $\lambda^*$ in \eqref{kkt}, and $z_i$ in \eqref{dy1} is given for the consensus of $\lambda_i$ to reach the same $\lambda^*$.  }

The communication noises $\epsilon_{ij}(k), ~\zeta_{ij}(k)$  can be used
to model information sharing uncertainties due to quantization errors (see \cite{peng}) or communication channel fading (see \cite{zhang2} and \cite{Zhang}).
Additionally, noises can be actively added  to achieve differential privacy protection as done  in \cite{DCU1}.

Define the $\sigma$-algebra at time $k$ as:
\begin{equation}\label{algebra}
\begin{split}
 \mathcal{F}_{k} =\sigma \{  &  \epsilon_{ij}(t), \zeta_{ij}(t),\delta_{i}(t),
\nu_i(t), L(t),  ~0 \leq t \leq k,  \\&~i,j=1,\cdots, N, ~ X(0), \Lambda(0), Z(0)\} .
\end{split}
\end{equation}
Define
$\mathcal{F}_{k}'=\sigma \{ \mathcal{F}_{k}, L(k+1 )\}$.
The following  assumptions imposed on  $\epsilon_{ij}(k), \zeta_{ij}(k), \delta_i(k), \nu_i(k)$ were also adopted in the existing SA and distributed optimization works (see \cite{Ned3}\cite{asum}\cite{ned}\cite{Zhang}).

\begin{assum}\label{assd}
For any $i \in \mathcal{N},$  $\{ \delta_{i }(k)\}$ is an i.i.d.  sequence  with zero mean and  bounded second moments $\sigma_{i,\delta}^{2}=E [\| \delta_{i }(k)\|^2].$
\end{assum}

\begin{assum}\label{asscn}
(i)  The  communication noises have conditional  zero mean,  i.e., $E[\zeta_{ij}(k) |\mathcal{F}_{k-1}']=\mathbf{0}$ and
$E[\epsilon_{ij}(k) |\mathcal{F}_{k-1}']=\mathbf{0}$.

(ii) There is a uniform bound on  conditional variances of the communication noise , i.e., there exists a constant $\mu>0$ such that
for any $i,j \in \mathcal{N}$ and any $k \geq 0$,
$E[\| \zeta_{ij}(k) \|^2 |\mathcal{F}_{k-1}'] \leq \mu^2$ and $ E[\|\epsilon_{ij}(k) \|^2 |\mathcal{F}_{k-1}'] \leq \mu^2.$

(iii)There exists a positive constant $c$ such that for any $i \in \mathcal{N} $ and any $k \geq 0$,
$$  E[ \nu_{i}(k) | \mathcal{F}_{ k-1}  ]=0,  ~~E[  \| \nu_{i}(k) \|^2  |  \mathcal{F}_{k-1} ]
\leq c (1+\|x_{i }(k)\|^2).$$

(iv) For all $  i\in \mathcal{N}  $, the sequences $\{L(k)\}$ and $\{ \delta_{i}(k)\}$   are mutually independent.
 The sequences $\{L(k)\}$ and $\{ \delta_{i}(k)\}_{ i\in \mathcal{N} }$   are   independent of $\mathcal{F}_{k-1}.$
  \end{assum}

\section{Convergence Analysis}

In this section,  we employ the ODE method  for SA algorithm to give  the convergence analysis for  algorithm \eqref{dy1}.
It is shown  with the following outline.
Theorem \ref{thm1} shows   that  the equilibrium point of the underlying ODE contains the optimal solution to  problem \eqref{problem}, while
Lemma \ref{lemode} shows the convergence of the underlying ODE.   Then Lemma \ref{lemnoise}  investigates properties of  the  extended noise sequences,
and Lemma \ref{lembound} shows that  the iteration sequence generated by   \eqref{dy1} are bounded.
Finally, Theorem \ref{thmcov} shows that the  estimates   generated  by \eqref{dy1}  converge  to the optimal  resource allocation with probability one.

%Denote by $X(k)=col\{x_1(k),...,x_n(k)\}$, $\Lambda(k)=col\{\lambda_1(k),...,\lambda_n(k)\}$, $Z(k)=col\{z_1(k),...,z_n(k)\}$, $\nabla f(X(k))=col \{ \nabla f_1(x_1(k)), \cdots, \nabla f_n(x_n(k)) \}$,   $D=col\{d_1,...,d_n\}$,
%$\nu(k)=col\{ \nu_{1}(k), \cdots,\nu_{N}(k)\}$, $\delta(k) =col\{ \delta_{1}(k), \cdots, \delta_{N}(k)\}$,    $\zeta(k)=col\{ \zeta_{1}(k), \cdots, \zeta_{N}(k)\}$ with $\zeta_{i}(k)=\sum_{j=1}^N a_{ij} (k) \zeta_{ ij} (k)$, and
%$\epsilon(k)=col\{ \epsilon_{1}(k) \cdots,\epsilon_{N}(k)\}$ with  $\epsilon_{i}(k)=\sum_{j=1}^N a_{ij} (k)\epsilon_{ij}(k)$.

Set $\zeta_{i}(k)=\sum_{j=1}^n a_{ij} (k) \zeta_{ ij} (k)$ and $\epsilon_{i}(k)=\sum_{j=1}^n a_{ij} (k)\epsilon_{ij}(k)$, and
\begin{equation}
\begin{array}{ll} \nonumber
& X(k)=col\{x_1(k),\cdots,x_n(k)\},               \quad               \Lambda(k)=col\{\lambda_1(k),\cdots,\lambda_n(k)\}, \\
& Z(k)=col\{z_1(k),\cdots,z_n(k)\},                \quad\;             \delta(k) =col\{ \delta_{1}(k),\cdots, \delta_{n}(k)\}, \\
& \nu(k)=col\{ \nu_{1}(k),\cdots,\nu_{n}(k)\}, \quad               D=col\{d_1,\cdots,d_n\},          \\
& \zeta(k)=col\{ \zeta_{1}(k),\cdots, \zeta_{n}(k)\}, \quad      \epsilon(k)=col\{ \epsilon_{1}(k),\cdots,\epsilon_{n}(k)\},\\
& \nabla f(X(k))=col \{ \nabla f_1(x_1(k)),\cdots, \nabla f_n(x_n(k)) \}.
\end{array}
\end{equation}
Then the   recursive  algorithm  \eqref{dy1}  can be rewritten in the compact form as follows:
\begin{equation}\label{cra}
\begin{array}{l}   \displaystyle X(k+1)           = P_{\Omega}  \big(   X(k)      + \alpha_{k} \big (    - \nabla f (X(k)) +\Lambda(k) - \nu(k)\big)\big),\\
\displaystyle   \Lambda(k+1)    =     \Lambda(k) +    \alpha_{k} \big(  -(L(k)\otimes I_m)( \Lambda(k)+ Z(k) ) \\
\displaystyle     \qquad \qquad \quad +     D- X(k) +\delta(k) + \zeta(k)+\epsilon(k)   \big) ,\\
\displaystyle Z(k+1)            =     Z(k)      + \alpha_{k}  \big(   (L(k)\otimes I_m) \Lambda(k) - \zeta(k) \big),  \end{array}
\end{equation}
where
$\Omega=\prod_{i=1}^n \Omega_i$ denotes the Cartesian product of $\Omega_i$.

Denote  by $e_1(k)= \big( (\bar{L}-L(k) )  \otimes I_m \big)( \Lambda _k  +  Z_k) ,$ $e_2(k)= \zeta(k) + \delta(k)+\epsilon(k)$, and by
$e_3(k)=\big( (L(k) -\bar{L} ) \otimes I_m \big) \Lambda(k) -\zeta(k)$. We then have
% Then \eqref{cra} is written as:
\begin{equation}
\begin{array}{l}
 \displaystyle X(k+1)          =   P_{\Omega}  \big( X(k)       + \alpha_{k}(    - \nabla f (X(k)) +\Lambda(k) - \nu(k)) \big),\\
 \displaystyle \Lambda(k+1)    =    \Lambda (k)  + \alpha_{k} \big(   - (  \bar{L}  \otimes I_m) ( \Lambda(k) +Z(k) )  \\
 \displaystyle  \qquad \qquad  +   D - X(k)+ e_1(k)+e_2(k) \big), \\
 \displaystyle Z(k+1)          =    Z(k)      + \alpha_{k}(  (  \bar{L}  \otimes I_m)  \Lambda(k) + e_3(k)). \end{array}
\end{equation}

By setting $S(k)= col\{ X(k),\Lambda(k), Z(k)\}$, we can regard the algorithm \eqref{cra} as an SA algorithm with  the following form:
\begin{equation}\label{compact}
S(k+1)  = P_{\Phi} \big( S(k) + \alpha_k (J(S(k)) +  \xi(k))\big),
\end{equation}
where
\begin{equation}\label{rootf}
\begin{split}
 & J(S) = \begin{pmatrix}
      &       -\nabla f(X) + \Lambda \\
      &     - (  \bar{L}  \otimes I_m) (  \Lambda + Z) + D-X\\
      &      (  \bar{L}  \otimes I_m)  \Lambda
\end{pmatrix},  \\& \xi(k)= \begin{pmatrix}
      &      - \nu(k)  \\
      &      e_1(k)+e_2(k) \\
      &   e_3(k)
\end{pmatrix}, \quad \Phi=\Omega \times \mathbf{R}^{mn} \times \mathbf{R}^{mn}.
\end{split}
\end{equation}

{
The convergence proof of \eqref{dy1} relies on the ODE  method for  SA (referring to \cite{sa} and \cite{sa2}).
Define the following  continuous-time projected dynamics as the underlying ODE of \eqref{dy1}
\begin{equation}\label{dy2}
 \dot{S}=J(S)+z, S(0)=col\{X(0),\Lambda(0),Z(0)\},
\end{equation}
with   $z\in -N_{\Phi} (S)$ being   the minimum force to keep the solution of \eqref{dy2} in $\Phi$, and $J(S)$ is defined by  \eqref{rootf}.

\begin{thm}\label{thm1}
Under Assumptions \ref{assp},\ref{assset}, and  \ref{assg}, \eqref{dy2} has at least one equilibrium point.
Furthermore, suppose $S^*=col \{X^*,\Lambda^*,Z^*\}$ is an equilibrium point of \eqref{dy2},  then   $S^*$   has  $X^*$ as   the optimal solution to problem  (\ref{problem}).
\end{thm}

{\bf Proof}:
Because problem \eqref{problem} is assumed to be solvable, there exist optimal solution $X^*$ and $\lambda^*\in \mathbf{R}^m$ such that \eqref{kkt} can be satisfied.
Then take $\Lambda=1_n\otimes \lambda^*$, $\bar{L}  \otimes I_m  \Lambda^*=\mathbf{0}$.
By $(\mathbf{1}^T_{n} \otimes I_m)  X^*= (\mathbf{1}^T_{n}\otimes I_m) D$
(that is $\sum_{i\in \mathcal{N}} x^*_i = \sum_{i\in\mathcal{ N}} d_i$), we have $D-X \in ker\{\mathbf{1}^T_{n}\otimes I_m\} $. Notice that $ker\{\mathbf{1}^T_{n}\otimes I_m\}$ and $ range\{\mathbf{1}_{n}\otimes I_m\} $    form
an orthogonal decomposition of $R^{nm}$ by the fundamental theorem of linear algebra.
Combined with  $ker(\bar{L} \otimes I_m)=range\{\mathbf{1}_{n}\otimes I_m\}$ due to Assumption \ref{assg}, we have   $D-X \in ker(\bar{L} \otimes I_m)^{\perp}$.
 Therefore, $D-X \in  range(\bar{L}  \otimes I_m )$, that is  there exists  $Z^*$ such that $ -\bar{L}  \otimes I_m Z^* + D-X^*=\mathbf{0}$. Hence, combined with \eqref{kkt},   $S^*=col\{X^*,\Lambda^*,Z^*\}$ is an equilibrium point of \eqref{dy2}.

On the other hand, when $S^*=col\{X^*,\Lambda^*,Z^*\}$ is an equilibrium point of \eqref{dy2}, it satisfies:
\begin{equation}
\begin{array}{lll}\label{eq2}
&& -\nabla f_i(x_i^*)+ \lambda_i^*  \in N_{\Omega_i}(x_i^*), x^*_i \in \Omega_i \\
&& (\bar{L}  \otimes I_m) (  \Lambda^* + Z^*) -( D-X^*)= \mathbf{0}\\
&& (\bar{L}  \otimes I_m)\Lambda^* = \mathbf{0}
\end{array}
\end{equation}

Since  $\bar{L}$  is the weighted Laplacian of an undirected  connected graph by Assumption \ref{assg}, it follows from $ (  \bar{L}  \otimes I_m)  \Lambda^*=\mathbf{0}_{mn}$ that
$\Lambda^*=  \mathbf{1}_n  \otimes \lambda^* $ for some $ \lambda^*\in \mathbf{R}^m$.
As a result, $\mathbf{0}_{m} \in   \nabla f_i(x_i^*)  -  \lambda^*  +  N_{\Omega_i}(x_i^*)$.
Furthermore,  $(\bar{L} \otimes I_m) \Lambda^*+(\bar{L} \otimes I_m)Z^* - (D-X^{*})=\mathbf{0}_{mn}$
 implies that  $(\bar{L} \otimes I_m)Z^*=D-X^*$. Then by noticing  $\mathbf{1}_n^T \bar{L}=\mathbf{0}_n^T$  we derive
 $\sum_{i\in \mathcal{N}} d_i=\sum_{i\in \mathcal{N}} x_i^*$. Moreover, $x_i^*\in \Omega_i$ due to the viability of ODE \eqref{dy2}.

Thus,   any equilibrium point  $S^{*}$ of \eqref{dy2} satisfies the KKT condition \eqref{kkt}, and
hence, $X^*$ is the optimal solution  to problem \eqref{problem}.
\hfill $\blacksquare$

Lemma \ref{lemode} shows that \eqref{dy2} converges to its equilibrium point  $S^*$.
\begin{lem}\label{lemode}
Under Assumptions  \ref{assp}, \ref{assset} and \ref{assg}, the trajectories of \eqref{dy2} are bounded
and  converge to its equilibrium point for any finite initial points.
\end{lem}

{\bf Proof}:
Take a Lyapunov function $V(S)=\frac{1}{2}||S-S^*  ||$, where $S^*$ is an equilibrium point of \eqref{dy2}.
Take $n_{\Omega}(X^*)\in N_{\Omega}(X^*)$ such that $\nabla f(X^*)-\Lambda^*+n_{\Omega}(X^*)=\mathbf{0}$, then
\begin{equation}
\begin{array}{ll}\label{eq3}
&\frac{dV}{dt}=(S -S^{*})^T (J(S) + z)\\
&\leq (X-X^*)^T(-\nabla f(X)+\Lambda +
\nabla f(X^*)-\Lambda^*+n_{\Omega}(X^*) ) \\
& + (\Lambda-\Lambda^*)^T(-\bar{L}\otimes I_m (\Lambda+Z)+(D-X)\\
&+\bar{L}\otimes I_m (\Lambda^*+Z^*)-(D-X^*)  ) + (Z-Z^*)^T \bar{L}\otimes I_m (\Lambda-\Lambda^*)\\
& \leq   -(X-X^*)^T(\nabla f(X)-\nabla f(X^*)) + (X-X^*)^Tn_{\Omega}(X^*)\\
 &-(\Lambda-\Lambda^*)^T \bar{L} (\Lambda-\Lambda^*)\leq 0
\end{array}
\end{equation}
Hence, any equilibrium point of \eqref{dy2} is Lyapunov stable, and given finite initial point $S(0)$,  the trajectories of \eqref{dy2} are bounded and belong to the compact forward invariant set $I_s=\{S | V(S)\leq V(0)\}$.

Denote $E$ as the set within $I_s$ such that $\dot{V}=0$.  Then we can show that the maximal invariance set in $E$ can only be $\{ S| \dot{S}=0\}.$
With the strict convexity of $f_i$,  $X=X^*$ must hold within set $E$. Furthermore, $\Lambda-\Lambda^* \in ker\{\bar{L}\otimes I_n\}$ by \eqref{eq3} and Assumption \ref{assg}.
Therefore, $\dot{Z}=\mathbf{0}_n$ and $Z=Z^*$ within set $E$. Moreover,  $\dot{\Lambda}=-L\otimes I_nZ^*+D-X^*$, and $\dot{\Lambda}$ must be $\mathbf{0}$; otherwise $\Lambda$ will go to infinity, which contradicts the boundedness of the trajectories.
Hence,  $\Lambda^* = \mathbf{1}_n \otimes \lambda^*$.
Therefore, all the trajectories of (\ref{dy1}) converge to the points in the maximal invariance set $\{ S| \dot{S}=0\}$.
Recalling the Lyapunov stability of $S^*$ and the LaSalle invariance principle, the dynamics (\ref{dy2}) converges to its equilibrium point $S^*$, which leads to the conclusion.
\hfill $\blacksquare$}

\subsection{ Extended noise property}

By  definition of $\mathcal{F}_k$ given in  \eqref{algebra},
$S(k)$  is adapted to $\mathcal{F}_{k-1 }$ according to \eqref{dy1}.
The extended noise sequence  $\{\xi(k)\}$ is state-dependent,
and  its properties  are shown in Lemma \ref{lemnoise}.

% About the noise sequence
\begin{lem}\label{lemnoise}
Suppose  Assumptions \ref{assg}, \ref{assd} and \ref{asscn} hold. Then
\begin{equation}\label{compon3}
E [\xi(k) | \mathcal{F}_{k-1}]=0,~~ E [ \| \xi(k) \|^2 |\mathcal{F}_{k-1}] \leq c_1 \|  S(k)\|^2+c_2~~a.s.
\end{equation}
for some finite constants $c_1,c_2$.
\end{lem}

{\bf Proof}:
By Assumption \ref{asscn} (iii),
\begin{equation}\label{conde0}
\begin{array}{lll}
 E[ \nu(k)| \mathcal{F}_{k-1}]          & =    & 0, \\
 E[ \| \nu(k) \|^2| \mathcal{F}_{k-1}]  & =    & \sum_{i=1}^n  E[ \| \nu_i(k) \|^2| \mathcal{F}_{k-1}]\\
                                              & \leq & nc+ c\| X(k)\|^2.
\end{array}
\end{equation}

Since  $a_{ij}(k)$ is adapted to  $\mathcal{F}^{'}_{k-1}$,  by   Assumption \ref{asscn} (i) we obtain \begin{equation}
 \begin{array}{lll}
 E[ \zeta_i(k)| \mathcal{F}_{k-1}']   & = & \sum_{j=1}^n  E [ a_{ij} (k) \zeta_{ ij} (k) |\mathcal{F}_{k-1}'] \\
                                    & = &  a_{ij} (k)  \sum_{j=1}^n  E [ \zeta_{ ij} (k) |\mathcal{F}_{k-1}'] =0, \nonumber
\end{array}
\end{equation}
and hence   $ E[ \zeta (k)| \mathcal{F}_{k-1}'] =0$.  By noting that   $\mathcal{F}_k \subset\mathcal{F}_k'$ we derive
$$ E[ \zeta (k)| \mathcal{F}_{k-1} ] =  E \big [ E[ \zeta (k)| \mathcal{F}_{k-1}']  \big | \mathcal{F}_{k-1} \big]=0. $$
Similarly,   it follows  that  $ E[ \epsilon (k)| \mathcal{F}'_{k-1}] =0$ and  $ E[ \epsilon (k)| \mathcal{F}_{k-1} ] =0$.

Since  $L(k)$  is independent of $\delta_i(k)$ and $\mathcal{F}_{k-1}$ by Assumption  \ref{asscn} (iv),  we obtain
$ E[  \delta_i (k)| \mathcal{F}_{k-1}' ] =E[  \delta_i (k)| \mathcal{F}_{k-1},L(k) ]=E[  \delta_i (k)| \mathcal{F}_{k-1} ].  $
Then by   Assumptions \ref{assd} and 4 (iv), we have  that, for any $i \in \mathcal{N}$
$E [\delta_i(k)|\mathcal{F}_{k-1}'  ]=E[  \delta_i (k)| \mathcal{F}_{k-1}]= E[  \delta_i (k)] =0.$
 Thus,  \begin{equation}\label{conditional1}
E [ e_2(k)| \mathcal{F}'_{k-1} ] = E[ \zeta (k)| \mathcal{F}'_{k-1} ] + E [\delta(k)|\mathcal{F}'_{k-1}  ] + E[ \epsilon (k)| \mathcal{F}'_{k-1} ] =0,
\end{equation}
which implies that $ E [ e_2(k)| \mathcal{F}_{k-1} ] =E[E[e_2(k)| \mathcal{F}'_{k-1}]\big | \mathcal{F}_{k-1}]=0.$

Note that $\Lambda(k)$ and $Z(k)$ are adapted to $\mathcal{F}_{k-1}$,  while
 $L(k)$ is independent of $\mathcal{F}_{k-1}$.  Then, by Assumption \ref{assg}
\begin{equation}
\begin{array}{lll}
 E[ e_1(k) | \mathcal{F}_{k-1} ]  & = &   E [ \big(  (\bar{L}-L(k) ) \otimes I_m \big) ( \Lambda (k)  +  Z(k)) |\mathcal{F}_{k-1} ]  \\
                                & = &   E [ \big(  (\bar{L}-L(k) ) \otimes I_m \big) |\mathcal{F}_{k-1} ] ( \Lambda (k)  +  Z(k)) \\
                                & = &   E [ \big(  (\bar{L}-L(k) ) \otimes I_m \big) ] ( \Lambda (k)  +  Z(k))=0, \\
 E[ e_3(k) | \mathcal{F}_{k-1} ]  & = &   E [   \big(  (\bar{L}-L(k) ) \otimes I_m \big) \Lambda(k) |\mathcal{F}_{k-1} ]  \\
                                & = &   E [  \big(  (\bar{L}-L(k) ) \otimes I_m \big)]    \Lambda(k)=0.
\end{array}
\end{equation}
Consequently, we conclude that $E [\xi(k) | \mathcal{F}_{k-1}]=0.$
\vskip 5mm
%and  $\mathcal{F}_k \subset\mathcal{F}_k'$ from,
Since $e_1(k)$ is adapted  to  $ \mathcal{F}_{k-1}'$ and $\mathcal{F}_k \subset\mathcal{F}_k'$ ,  it follows from \eqref{conditional1} that
\begin{equation}\label{conde01}
\begin{split}
 E[ e_1(k)^T e_2(k) | \mathcal{F}_{k-1} ]  &=   E \big [ E[ e_1(k)^T e_2(k) | \mathcal{F}_{k-1}'] \big | \mathcal{F}_{k-1}  \big]   \\
                                         & =   E \big [ e_1(k)^T  E [e_2(k) | \mathcal{F}_{k-1}^{'} ] \big | \mathcal{F}_{k-1} \big]  =0. \end{split}
\end{equation}

Since $\Lambda(k)$  and $L(k)$ are adapted to $\mathcal{F}_{k-1}'$, by $ E[ \zeta (k)| \mathcal{F}_{k-1}'] =0$  and  $\mathcal{F}_k \subset\mathcal{F}_k'$, we get
\begin{equation}\label{conde1}
\begin{split}
&E [   \big(  (L(k) -  \bar{L} )  \otimes I_m \big)  \Lambda(k) \big)^T  \zeta(k) |\mathcal{F}_{k-1} ]  \\
&\quad = E \Big[ E \big [ \big(  (L(k) -  \bar{L} )  \otimes I_m \big)  \Lambda(k) \big)^T  \zeta(k) | \mathcal{F}_{k-1}^{'}\big]  \Big|\mathcal{F}_{k-1} \Big ]  \\
&\quad =  E \Big[  \big(  (L(k) -  \bar{L} )  \otimes I_m \big)  \Lambda(k) \big)^T   E \big [  \zeta(k) | \mathcal{F}_{k-1}^{'}\big]  \Big|\mathcal{F}_{k-1} \Big ] =0.
\end{split}
\end{equation}

By the conditional H\"older   inequality
 \begin{equation}\label{Holder}
E[\| X^TY\|   \big | \mathcal{F}]  \leq(E[\| X \|^2  \big | \mathcal{F}] )^{\frac{1}{2}} (E[\| Y \|^2  \big |  \mathcal{F}] )^{\frac{1}{2}}  ,
\end{equation}
from  Assumption  \ref{asscn}  (ii)  we see  that
 \begin{equation}
 \begin{array}{lll}
 & E[ \zeta_{ij}(k)  ^T \zeta_{ip}(k)    \big |  \mathcal{F}_{k-1}'  ] \leq
E[ \| \zeta_{ij}(k)  ^T \zeta_{ip}(k) \|  \big |  \mathcal{F}_{k-1}'  ]  \\&\leq(E[\| \zeta_{ij}(k)  \|^2  \big |  \mathcal{F}_{k-1}'] )^{\frac{1}{2}}
(E[\| \zeta_{ip}(k)  \|^2   \big |  \mathcal{F}_{k-1}'] )^{\frac{1}{2}}  \leq \mu^2. \nonumber
 \end{array}
\end{equation}
Then, since $A(k)$ is   adapted to  $\mathcal{F}_{k-1}'$, we have
\begin{equation}
\begin{array}{lll}
E [\|\zeta_i(k) \|^2 |  \mathcal{F}_{k-1}'  ] & = & E[\sum_{j,p =1}^n a_{ij}(k) a_{ip}(k)\zeta_{ij}(k)^T \zeta_{ip}(k)  |  \mathcal{F}_{k-1}'  ] \\
                      & = & \sum_{j,p =1}^n a_{ij}(k) a_{ip}(k) E[ \zeta_{ij}(k)^T \zeta_{ip}(k) |  \mathcal{F}_{k-1}'   ] \\
                      & \leq &\sum_{j,p =1}^n \mu^2=n^2\mu^2 . \nonumber
\end{array}
\end{equation}
Similarly, $
E[\| \epsilon_i(k) \|^2  |  \mathcal{F}_{k-1}'  ] \leq n^2 \mu^2~ \forall i \in \mathcal{N}.$
From  Assumption  \ref{asscn} (iv), it is clear that $\delta_{i}(k)$ is independent of $\mathcal{F}_{k-1}'$,
and hence, by Assumption \ref{assd}
  $E [\| \delta(k) \|^2 | \mathcal{F}_{k-1}' ]=E [\| \delta(k) \|^2]=\sum_{i=1}^n E [\| \delta_i(k) \|^2  ] .$
In summary,
\begin{equation}\label{sndc}
\begin{array}{ll}
&E [\| \zeta(k) \|^2   |  \mathcal{F}_{k-1}'  ]     \leq  n^3 \mu^2, ~~
E [\| \epsilon(k) \|^2 |  \mathcal{F}_{k-1}'   ] \leq n^3\mu^2, \\
&E [\| \delta(k) \|^2   |   \mathcal{F}_{k-1}'   ]    =  \sum_{i=1}^n\sigma_{i,\delta} \triangleq\sigma_{\delta} .
\end{array}
\end{equation}
Then
\begin{equation}
 \begin{array}{lll}
& E [ \|  e_2(k)\|^2 |\mathcal{F}_{k-1}']   \\ & \leq  3\big(   E [\| \epsilon(k) \|^2 |\mathcal{F}_{k-1}'] +E [\| \zeta(k) \|^2 |\mathcal{F}_{k-1}'] +E [\| \delta(k) \|^2 |\mathcal{F}_{k-1}'] \big) \\ & \leq 3(2n^3\mu^2+\sigma_{\delta }^2)  \triangleq C_{1}, \nonumber
 \end{array}
\end{equation}
and hence, by  $\mathcal{F}_k \subset\mathcal{F}_k'$, we get
\begin{equation}\label{conde2}
E [ \|  e_2(k)\|^2 |\mathcal{F}_{k-1}]=E \big[E [ \|  e_2(k)\|^2 |\mathcal{F}_{k-1}']   \big | \mathcal{F}_{k-1}\big]  \leq C_{1}.
\end{equation}

Because  $ \Lambda (k)  $ and $ Z(k)$ are adapted to $\mathcal{F}_{k-1}$,   and  $L(k)$  is independent of $\mathcal{F}_{k-1}$ by Assumption 4 (iv), we have
\begin{equation}
\begin{array}{lll}
E [ \| e_1(k) \|^2 |\mathcal{F}_{k-1}]  &=E [ \|   \big( (\bar{L}-L(k) )  \otimes I_m \big) ( \Lambda (k)  +  Z(k))   \|^2 |\mathcal{F}_{k-1}]  \\
& \leq C_{2 } \| \Lambda(k) + Z(k)\|^2, \nonumber
\end{array}
\end{equation}
where $C_2=E [\|  L(k) - \bar{L} \|^2]$  is finite.
It, along with \eqref{conde01} \eqref{conde2}, yields
\begin{equation}\label{compon1}
\begin{array}{lll}
&&E [ \| e_1(k) +e_2(k)\|^2 |\mathcal{F}_{k-1}]  =  E [ \|  e_1(k) \|^2 |\mathcal{F}_{k-1}] \\
&&  + E [ \|  e_2(k)\|^2 |\mathcal{F}_{k-1}] + 2E [   e_1(k)^Te_2(k)   |\mathcal{F}_{k-1}]   \\
&&\leq C_{2 } \| \Lambda(k) + Z(k)\|^2 +C_1 .
\end{array}
\end{equation}

By \eqref{conde1} \eqref{sndc} and  $\mathcal{F}_k \subset\mathcal{F}_k'$, we derive
\begin{equation}\label{compon2}
\begin{array}{lll}
E [ \| e_3(k) \|^2 |\mathcal{F}_{k-1}] & =  & E [ \|  \big( (L(k) -\bar{L} )  \otimes \big) \Lambda(k) -\zeta(k)  \|^2 |\mathcal{F}_{k-1}]  \\
                                     & =  & E [ \| \big(  (L(k) -\bar{L} )  \otimes I_m\big) \Lambda(k) \|^2 |\mathcal{F}_{k-1}]  \\
                                     & \; & - 2E [   \big( \big( (L(k) -\bar{L} )  \otimes I_m \big) \Lambda(k) \big)^T\zeta(k)    |\mathcal{F}_{k-1}]   \\
                                     & \; & + E [ \|  \zeta(k)  \|^2 |\mathcal{F}_{k-1}]   \leq C_2 ||\Lambda(k) ||^2+n^3\mu^2.
\end{array}
\end{equation}

 In summary, from \eqref{conde0},  \eqref{compon1}, and  \eqref{compon2}, we obtain
 \begin{equation}\label{compon3}
\begin{array}{lll}
  E [ \| \xi(k) \|^2 |\mathcal{F}_{k-1}] & =    &  E [ \|  \nu(k) \|^2 |\mathcal{F}_{k-1}]+E [ \| e_3(k) \|^2 |\mathcal{F}_{k-1}]\\
                                       & \;   & +E [ \| e_1(k) +e_2(k)\|^2 |\mathcal{F}_{k-1}]\\
                                       & \leq & nc+ c\| X(k)\|^2 +   C_2 ||\Lambda(k) ||^2  + n^3\mu^2\\
                                       &  \;  &+ C_2\|\Lambda(k) + Z(k)\|^2 +C_1  \\
                                       & \leq & c_1 \|  S(k)\|^2+c_2 \nonumber
\end{array}
\end{equation}
for some positive constants $c_1, c_2$.
\hfill $\blacksquare$

\subsection{Stability}
The following result is about the boundedness of the iterations  before showing its convergence.

\begin{lem}\label{lembound}
Under Assumptions 1-4, $\{ S(k)\}$ generated by the distributed algorithm \eqref{dy1}  is bounded with probability one given any finite initial value $S(0)$.
\end{lem}

{\bf Proof:}
Denote by $S^{*}$ as an equilibrium point of \eqref{dy2}, i.e., $J(S^{*})\in N_{\Phi}(S^*)$.
Then, by Assumption  \ref{assp} and  the KKT condition \eqref{kkt},  $S^{*}$ is a  finite value.
Take $v(S)= \| S- S^{*}\|^2$ as a Lyapunov function. {  Then from  \eqref{compact} and  the non-expansive property of the projection operator \eqref{pro} we derive}
\begin{equation}
\begin{array}{ll}
v( S(k+1)) & =  \| S(k+1)- S^{*}\|^2\\
           & \leq  \| S(k) -S^{*}+ \alpha_k (J(S(k)) +  \xi(k))\|^2 \\
           & \leq  \| S(k )- S^{*}\|^2+ 2\alpha_k (S(k) -S^{*})^T (J(S(k)) +  \xi(k)) \\ &+ \alpha_k^2\big(\|   J(S(k))\|^2 +2 \xi(k)^TJ(S(k)) +\|  \xi(k) \|^2\big) . \nonumber
\end{array}
\end{equation}
 Since $S(k)$ is adapted to  $\mathcal{F}_{k-1 }$,    by   recalling  $E [\xi(k) | \mathcal{F}_{k -1}]=0$ from Lemma \ref{lemnoise} we obtain
 \begin{equation}\label{inequality0}
\begin{array}{lll}
 E [ v( S(k+1)) | \mathcal{F}_{k-1}]  & \leq & v(S(k))+ 2\alpha_k (S(k) -S^{*})^T J(S(k)) \\
                                  & +    & \alpha_k^2 ( \|   J(S(k)) \|^2  + E [  \|  \xi(k) \|^2 | \mathcal{F}_{k-1}]). \nonumber
\end{array}
\end{equation}
Similar to the proof of Lemma \eqref{lemode}, $(S(k) -S^{*})^T J(S(k))  \leq    0.$
Then  by  Lemma  \ref{lemnoise}, we get
\begin{equation}\label{ineq1}
E [ v( S(k+1)) | \mathcal{F}_{k-1}]  \leq v(S(k)) + \alpha_k^2 ( \|   J(S(k)) \|^2  + c_1 \|  S(k)\|^2+c_2) .
\end{equation}
 {  From  Assumption \ref{assp} and taking $n_{\Omega}(X^*)\in N_{\Omega}(X^*)$ such that $\nabla f(X^*)-\Lambda^*+n_{\Omega}(X^*)=\mathbf{0}$, we have
\begin{equation} \label{inequality01}
\begin{array}{ll}
 \|   J(S(k)) \|^2  & =   \| -\nabla f(X(k))+\Lambda(k) +\nabla f(X^{*}) -\Lambda^{*} +n_{\Omega}(X^*)\|^2 \\
                    &  +\|  (  \bar{L}  \otimes I_m)  \big( (Z(k)-Z^{*})+ (\Lambda(k)-\Lambda^{*}) \big )\\
                    & +    X(k)-X^{*}      \|^2  +  \|    (  \bar{L}  \otimes I_m)   (\Lambda(k)-\Lambda^{*}) \|^2 \\
                    &\leq 3(\| \nabla f(X(k)) -\nabla f(X^{*})  \|^2 + \|  \Lambda(k)- \Lambda^{*}\|^2\\
                    &+ \| n_{\Omega}(X^*)\|^2+ \|   (  \bar{L}  \otimes I_m)  (\Lambda(k)-\Lambda^{*}) \|^2 +\|    X(k)-X^{*}   \|^2 \\
                    & + \|  (  \bar{L}  \otimes I_m)  (Z(k)-Z^{*})  \|^2 )+  \|    (  \bar{L}  \otimes I_m)   (\Lambda(k)-\Lambda^{*}) \|^2\\
                    &\leq (3l^2_c+3)  \|    X(k)-X^{*} \|^2+3c_4\| Z(k)-Z^{*}\|^2\\
                    &+(3+4c_4) \| (\Lambda(k)-\Lambda^{*}) \|^2 +c_n \\
                    &\leq (3+ 3l_c^2+4c_4 ) ||S(k)-S^{*}||^2+c_n = c_5 v( S(k))+c_n,
\end{array}
\end{equation}
where  $c_4=\| \bar{L}\|$ and $c_n= \| \nabla f(X^*)-\Lambda^*\|^2$.}
Note that  $$\|  S(k)\|^2 \leq 2(\|  S(k)-S^{*}\|^2 +\|  S^{*}\|^2)=2(v( S(k))\|^2 +\|  S^{*}\|^2).$$    Incorporated with \eqref{ineq1} and \eqref{inequality01},  it  yields
\begin{equation}
 \begin{array}{lll}
 & E [ v( S(k+1)) | \mathcal{F}_{k-1}]  \leq  v(S(k)) \\
                                  &  +             \alpha_k^2 \big( c_5v( S(k))+c_n +2c_1v(S(k)) + 2c_1\|  S^{*}\|^2+c_2\big)\\
                                  & \triangleq    (1+ c_6 \alpha_k^2)   v(S(k))+ c_7\alpha_k^2 ,
\end{array}
\end{equation}
where $c_6=2c_1+c_5 , c_7=2c_1\|  S^{*}\|^2+c_2+c_n$.

Since  $\{\alpha_k\}$  satisfies \eqref{stepsize}, with probability one
 $\lim\limits_{k \rightarrow \infty} v(S(k))$  exists and is finite  by Lemma \ref{marg} in Appendix. Therefore,   $\{ S(k)\}$ is bounded with probability one.
\hfill $\blacksquare$

\subsection{Convergence}
The following result gives the main convergence  result  for the SA-based distributed algorithm \eqref{dy1}.

\begin{thm}\label{thmcov}
Suppose   Assumptions \eqref{assp}-\eqref{asscn} hold.
Let sequences  $ \{ x_{i}(k)\},~ \{ \lambda_{i}(k) \},~ \{ z_{i}(k) \}$ be produced by \eqref{dy1}
given  any finite initial values $x_{i}(0),~\lambda_{i}(0),~ z_{i}(0)$. Then
$$ \lim_{k \rightarrow \infty} x_{i }(k)=x_i^{*}~~ a.s.,$$
  where $X^{*}=col\{x_1^{*}, \cdots, x_n^{*}\}$ is the optimal  resource allocation to problem \eqref{problem}.
\end{thm}

{\bf Proof:}
{  Note that $\theta_k$, $Y_k$, $g(\theta)$  and $\Phi$ for  \eqref{constrained} correspond  to
$\theta_k=S(k) $,  $Y_k=J(S(k)) +  \xi(k)$, $g(\theta)=J(S)$ and
$\Phi=\Omega \times \mathbf{R}^{mn} \times \mathbf{R}^{mn}$ for   \eqref{compact}.
Then we can  apply Theorem \ref{CSA}  in Appendix to prove the conclusion,  and  it suffices to check conditions C1-C4 given in Appendix.

Since $S(k)$ is adapted to $\mathcal{F}_{k-1}$, by  \eqref{compon3}   we drive
$$E[\|Y_k\|^2|\mathcal{F}_{k-1}] \leq  \| J(S(k)) \|^2+c_1 \|  S(k)\|^2+c_2~~a.s.$$
Then by  Lemma \ref{lembound}, \eqref{inequality01} and   Assumption \ref{assp}  we conclude that C1 hold.
From  \eqref{rootf} and Lemma \eqref{lemnoise} it is easily seen that C2 holds.
By definition of  $J(S)$  given by \eqref{rootf} and Assumption \ref{assp} we know that C3 holds.
Since  $\{ S(k)\}$  is bounded with probability one from  Lemma \ref{lembound}, we  then have C4.

As a result, C1-C4 hold.  Since $\Phi=\Omega \times \mathbf{R}^{mn} \times \mathbf{R}^{mn}$, with Assumption \ref{assset}  it is easily seen  that $\Phi$   satisfies the similar conditions as $\Omega_i.$
 Then, by  Theorem \ref{CSA},  $S(k)$ converge with probability one  to
  the invariant set of  \eqref{dy2}.
Thus, by Lemma \ref{lemode}, $X(k)$ converges with probability one to the optimal solution $X^{*}$.
}
 \hfill $\blacksquare$

\section{Demand Response Management and Simulations}

In this section, we apply the  RA optimization model  \eqref{problem} and algorithm \eqref{dy1}
 to distributed multi-period demand response  management in  power systems (see \cite{stevon2} and \cite{peng2}).

Suppose that a group of load aggregators (with index
$\mathcal{N}=\{1,\cdots,n\}$) need to decide the load demand in the following $T$
periods $P_i^d\in \mathbf{R}^T$, in order to meet the generation
scheduling $P_i^g \in \mathbf{R}^T$ and minimize the disutilities.
$P_i^g$ is usually decided  by other decision processes based on
the generator unit commitment or real-time generation prediction of
renewables, which is fixed and assumed  to  be only informed or
observable to  agent $i$.  Aggregator  $i$ formulates its local
objective function $f_i(P^d_i)$ to consider the  costs  or
disutilities  due to demand response $P_i^d$. { Moreover,
$P_i^d\in \Omega_i$ specifies the local response constraints, which
considers the lower and upper bounds in each period, the total
demand in the following  $T$ periods, ramping constraints, and
other local specifications.}  Hence,  the multi-period demand
response management problem is formulated as:
\begin{equation}\label{problem2}
\begin{array}{lll}
\min_{P^d_i \in \mathbf{R}^T,i\in \mathcal{N}} & \; & \sum_{i\in \mathcal{N}} f_i(P^d_i)\\
subject \; to \;     &\;& \sum_{i\in \mathcal{N}} P^d_i = \sum_{i\in\mathcal{ N}} P^g_i, \quad P_i^d \in \Omega_i
\end{array}
\end{equation}

In many practical cases,  $P_i^g$ can only be observed indirectly
through local measurements of wind speed, or solar radiation, or
local frequency deviation, and hence,  suffers from various
observation noises. In addition,  $f_i(P_i^d)$ should take full
consideration of user's demand requirements,  (dis)utility,
satisfactory levels, and payoffs, and hence,  is influenced by
various external factors, such as temperature, electricity price,
and renewable generations. Therefore, the gradient observation of
$f_i(P_i^d)$ may also be noisy. The aggregators may share
information through wireless communication networks with switching
topologies and noisy channels. As a result, algorithm \eqref{dy1}
can be applied to handle the above challenges for problem
\eqref{problem2}. Compared with previous works \cite{stevon2} and
\cite{peng2}, the proposed model  here considers the demand response
in multi-periods and local load response feasibility constraints,
and the algorithm can handle various observations and communication
uncertainties, which may be more practical in many cases.

In what follows, we give a numerical experiment to illustrate the
algorithm performance.

\begin{exa}
Consider the following three-period  demand response management problem:
\begin{equation}\label{exam_problem}
\begin{array}{lll}
&&\min_{P^d_i \in \mathbf{R}^{3}, i\in \mathcal{N} } \sum_{i\in \mathcal{N}}  E_{\Psi_i,\theta_i}[ {P_i^d}^T (Q_i + \Psi_i) P_i^d  +  (c_i + \theta_i)^T P^d_i] \\
&& s.t. \qquad \qquad      \sum_{i\in \mathcal{N}} P^d_i =\sum_{i \in \mathcal{N}} P^g_i \\
&&      \quad \qquad \qquad       R_i P_i^d \leq l_i, R_i \in \mathbf{R}^{12\times 3}, l_i\in \mathbf{R}^{12\times 1}, i\in \mathcal{N},
\end{array}
\end{equation}
{  where $R_i P_i^d \leq l_i$  is the compact form of  the following local load response feasibility constraints:
$ [l_i]_{11} \leq \mathbf{1}^T P_i^d \leq [l_i]_{21}$,
$[l_i]_{31}\leq [P_i^d]_{11}-[P^d_i]_{21}\leq [l_i]_{41} $,
 $[l_i]_{51}\leq [P_i^d]_{21}-[P_i]_{31}\leq [l_i]_{61}$,
 $ [l_i]_{71} \leq  [P_i^d]_{11} \leq [l_i]_{81} $,
 $ [l_i]_{91} \leq [P_i^d]_{21} \leq [l_i]_{10,1} $ and
 $ [l_i]_{11,1} \leq [P_i^d]_{31} \leq [l_i]_{12,1} $. }

The basic simulation experiment settings are given as follows. The number of  agents
is set to be $10$. $Q_i $ and $c_i$ are randomly generated symmetric
positive definite matrices and random  vectors, respectively. Each
$P^g_i$ and $l_i$ are also randomly generated vector that can ensure
Assumptions \ref{assp} and \ref{assset}.
\begin{figure}
  \centering
  \includegraphics[width=3in]{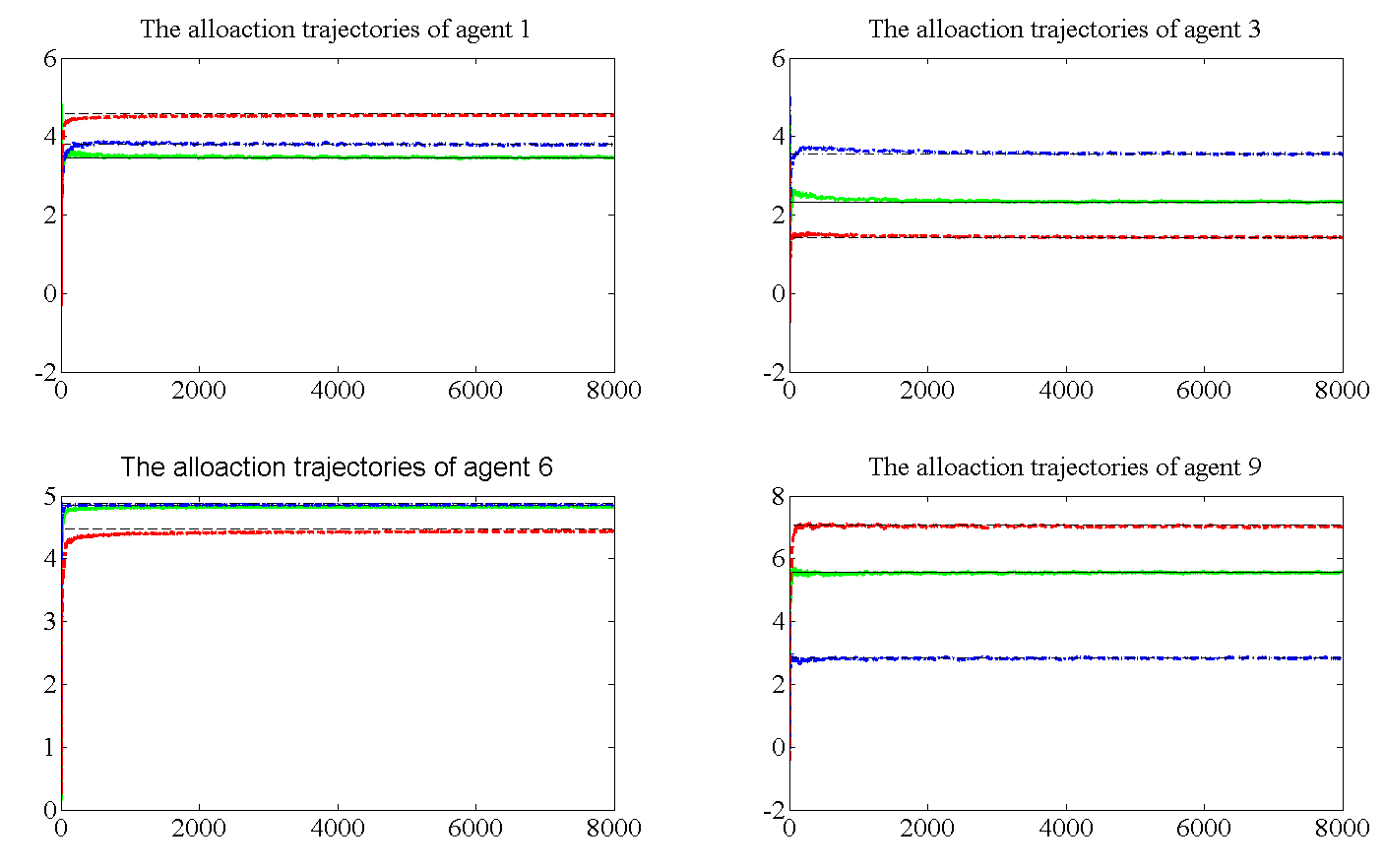}\\
  \caption{ The averaging trajectories of some agents' allocation variables
  }\label{fig_max1}
\end{figure}

Consider a graph set $\mathcal{G}_s$  containing  $30$ graphs, each
of which is generated  according to the random graph model
$G(10,P)$, where $P$ is the probability of occurrence for  any
possible edge. The probability $P$ is randomly and uniformly drawn
from $[0.05,0.1]$ for each graph in $\mathcal{G}_s$. Select a graph
set $\mathcal{G}_s$ with its union graph being connected. At time $k$,  a graph is randomly drawn from the graph
set $\mathcal{G}_s$ according to the uniform distribution.

For $i \in \mathcal{N}$, $[\Psi_i]_{ij},~[\theta_i]_j $ are i.i.d.
random variables satisfying the Gaussian distribution $N(0,0.5)$
with zero mean and variance $0.5$. Let both the generation
observation noise $\delta_i$ and communication noise $\zeta_{ij}$,
$\epsilon_{ij}$  be  i.i.d. random vectors satisfying the Gaussian
distribution $N(\mathbf{0},I_3)$ with zero mean vector and
covariance matrix $I_3$. Hence, Assumptions \ref{assd} and
\ref{asscn} are satisfied. The stepsize $\alpha_k$ in \eqref{dy1} is
set as $\alpha_k=\frac{1}{(k+1)^{0.6}}$.

{\bf Experiment 1}: Given a randomly generated graph set
$\mathcal{G}_s$ and a randomly generated setting for problem \eqref{exam_problem}, we apply
algorithm \eqref{dy1} to generate $200$ independent sample paths
with iteration length of $8000$.
Figure \ref{fig_max1} shows the averaging trajectories of some
agents' allocation variables, and illustrates how the agents find the
optimal allocation. Moreover, Figure \ref{fig_max2} shows the
averaging trajectories of some algorithm performance indexes,
including the distance to optimal solution $||P^d-{P^d}^* ||$,
function value $f(P^d)$, $||\bar{L}\Lambda ||$, and $ || \sum_{i\in
\mathcal{N}} (P_i^d-P_i^g)||$.

\begin{figure}
  \centering
  \includegraphics[width=3in]{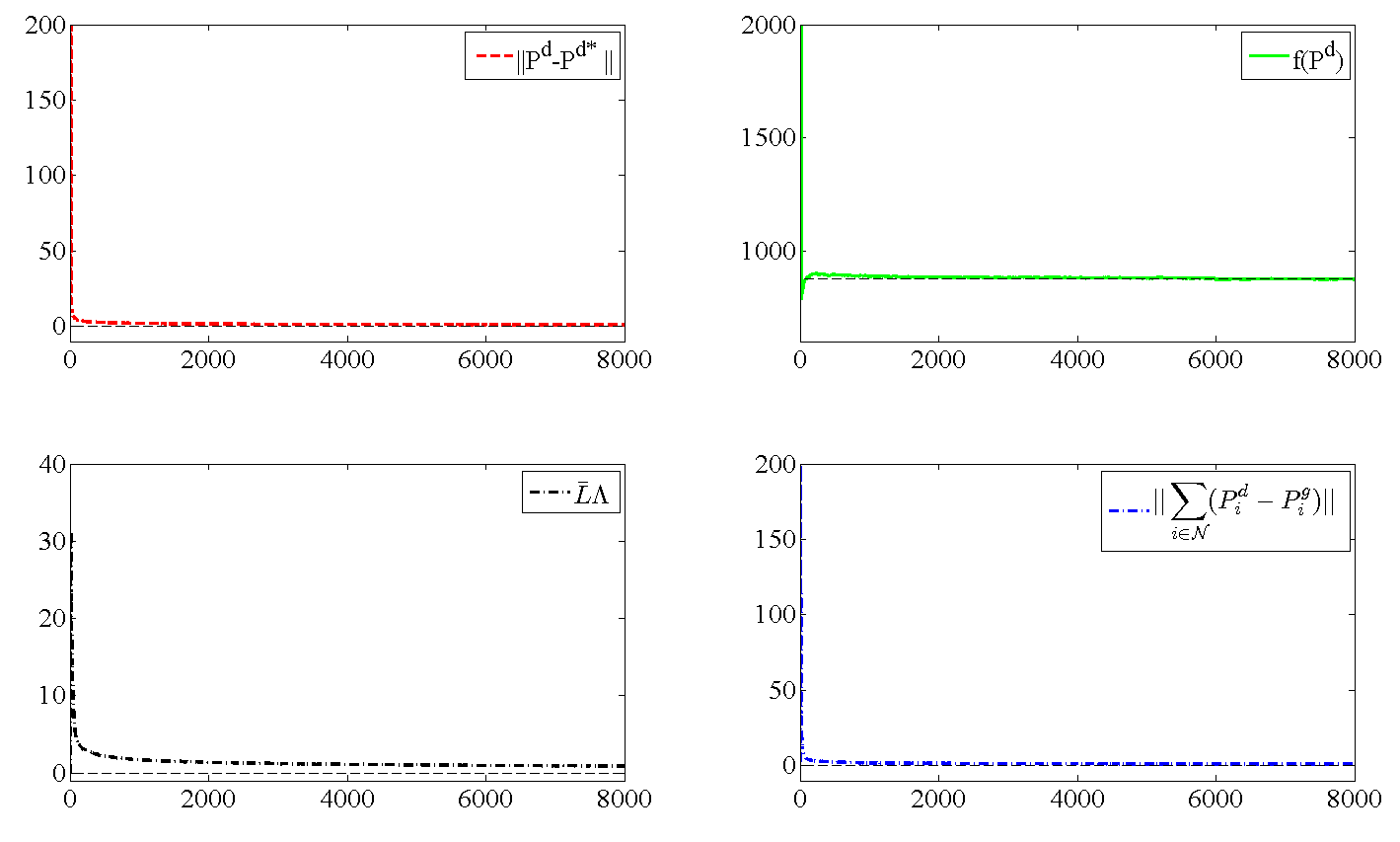}\\
  \caption{  The averaging trajectories of some performance indexes.
  }\label{fig_max2}
\end{figure}
{\bf Experiment 2}:  Let us randomly generate a graph set
$\mathcal{G}_s$  and a setting for problem \eqref{exam_problem} at each round of this simulation, and
employ algorithm \eqref{dy1} to generate one sample path of this setting with
iteration length of $8000$. We repeat the procedure for $100$ rounds,
and use Figure \ref{fig_max3} to show  the histograph of some
performance indexes at iteration time $8000$. It illustrates that
algorithm \eqref{dy1} can almost surely find the optimal allocation
for different problem  settings with only one sample path.

\begin{figure}
  \centering
  \includegraphics[width=3in]{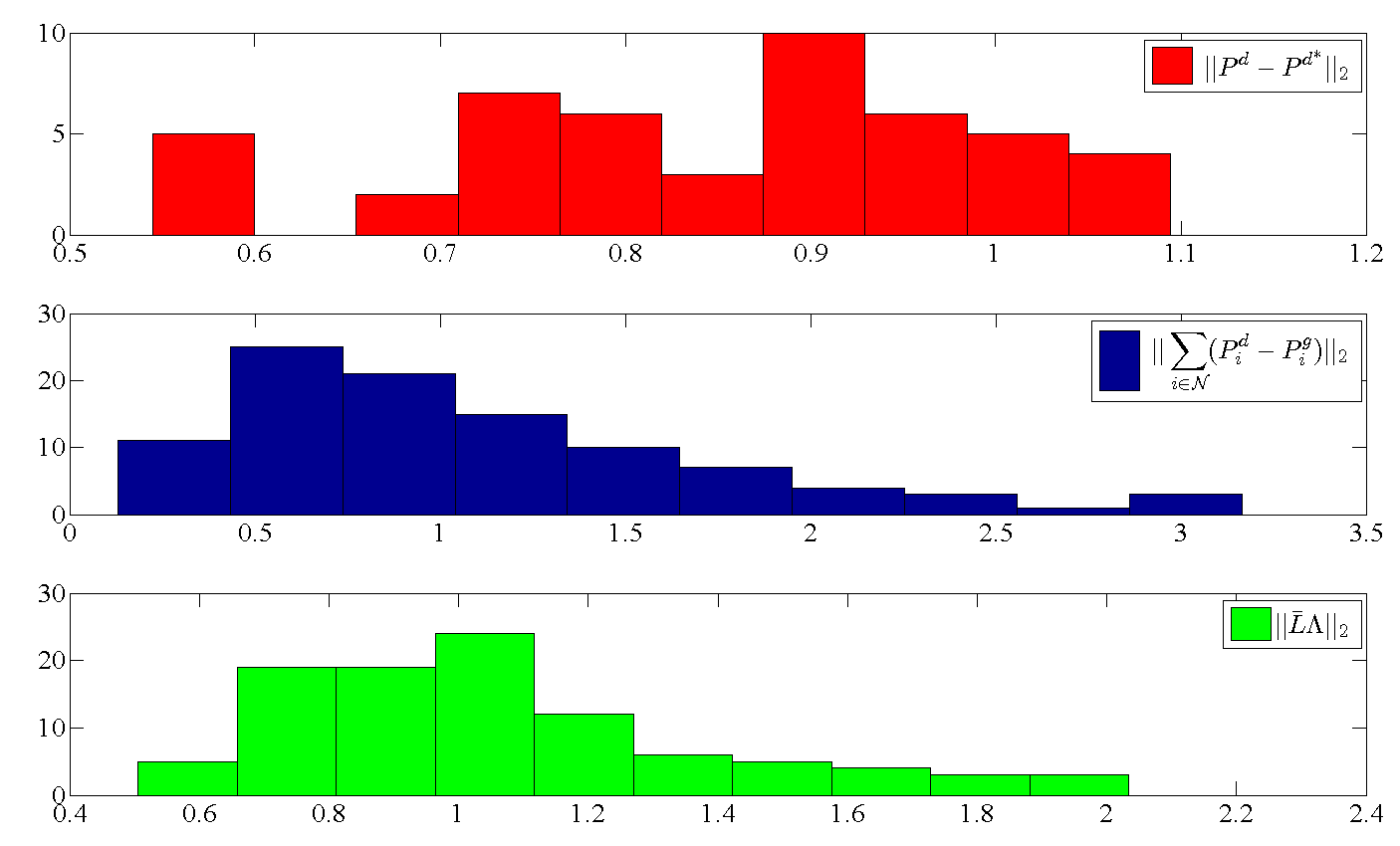}\\
  \caption{  The histograph of some performance indexes at iteration time $k=8000$.
  }\label{fig_max3}
\end{figure}

\end{exa}

\section{Conclusions}
In this paper, an SA-based distributed algorithm was proposed to solve
a class of RA optimization problems under various uncertainties.
The gradient and resource observation noises were taken into consideration, and
the communication network was assumed with randomly switching topologies and noisy communication channels.
The algorithm was proved to converge to the optimal solution with probability one by resorting to the ODE method for SA algorithm,
which may
%(The following sentence is newly added.)
demonstrate great potentials of SA algorithm and ODE methods for distributed decision problems over network systems under noisy data observations.

{
\section*{Appendix}%{Results on Constrained Stochastic Approximation and Supermartingale}
Here is the  convergence result for the constrained stochastic
approximation.
%The symbols used here are within the appendix (you have used C1-C4 out of this appendix).
Consider
\begin{equation}
\label{constrained}
\theta_{k+1}=P_{\Phi} \{ \theta_k+\alpha_kY_k\},
\end{equation}
where $\Phi \in \mathbf{R}^m$ is a convex constraint set.
Next follows the conditions  for its convergence analysis.

C1:  $\sup_k E [ \| Y_k\|^2]<\infty. $

C2: There is a measurable  function $g(\cdot)$ such that
$$ E_k[Y_k]=E[Y_k|\theta_0,Y_i, i<k]=g(\theta_k) .$$

C3: $g(\cdot)$ is continuous.

C4: $\theta_k$ is bounded  with probability one.

\begin{thm} \label{CSA}\cite[Theorem 5.2.1 and Theorem 5.2.3] {sa} Let C1-C4, and    \eqref{stepsize} hold for algorithm
\eqref{constrained}.     If  $\Phi$   satisfies  the same condition
as  that  imposed on $\Omega_i$ in Assumption \ref{assset}, then
with probability one   $\theta_k$ converges  to the invariant set of
the following projected ODE in $\Phi$:
$$\dot{\theta}=g(\theta)+z ,$$
where   $z\in -N_{\Phi} (\theta)$ is  the minimum force to keep the trajectories of the projected ODE in $\Phi.$

 \end{thm}
}

The following lemma shows convergence properties for nonnegative  super-martingales.
\begin{lem}[Robbins-Siegmund](\cite{supermds}) \label{marg}
Let $(\Omega,\mathcal{F},\mathbb{P})$ be a probability space and $\mathcal{F}_0\subset \mathcal{F}_1 \subset\cdots $
be a sequence of $\sigma-$algebra of $\mathcal{F}$. Let $\{d_k\}$ and
$\{w_k\}$ be nonnegative  $\mathcal{F}_{k}$-measurable  random variables such that
$$ E[d_{k+1}| \mathcal{F}_{k}] \leq (1+\alpha_k) d_k + w_k  ,$$
where $\alpha_k \geq 0$ are deterministic scalars with  $\sum_{k=1}^{\infty}\alpha_k <\infty$.
If  $\sum_{k=1}^{\infty} w_k < \infty$, then  $\{d_k\}$
converges  with probability one to some finite  random variable.
\end{lem}

\end{document}